\documentclass[11pt]{article}

\usepackage{latexsym}
\usepackage{amssymb}

\newtheorem{thm}{Theorem}

\newtheorem{rmr}{Remark}
\newtheorem{lem}{Lemma}

\begin{document}
{
\begin{center}
{\Large\bf
Invertible extensions of symmetric operators and the corresponding generalized resolvents.}
\end{center}
\begin{center}
{\bf S.M. Zagorodnyuk}
\end{center}

\section{Introduction.}
Let $A$ be a closed symmetric invertible operator in a Hilbert space $H$.
The domain of $A$ is not supposed to be necessarily dense in $H$.
Let $\widetilde A$ be a self-adjoint extension of $A$, acting
in a Hilbert space  $\widetilde H\supseteq H$.
Recall that an operator-valued function $\mathbf{R}_\lambda$, given by the following relation:
$$ \mathcal{\mathbf{R}}_\lambda = \mathbf{R}_\lambda(A) = \mathbf{R}_{s;\lambda}(A) =
P^{\widetilde{H}}_H \left( \widetilde A - \lambda E_{\widetilde{H}} \right)^{-1}|_H,\qquad
\lambda\in \mathbb{R}_e, $$
is said to be {\bf a generalized resolvent}
of the symmetric operator $A$ (corresponding to the extension $\widetilde A$).
Fix an arbitrary point $\lambda_0\in \mathbb{R}_e$.
An arbitrary generalized resolvent $\mathbf{R}_{s;\lambda}$ of the operator $A$
is given by Shtraus's formula:
\begin{equation}
\label{f11_3_p2_1}
\mathbf R_{s;\lambda} = \left\{ \begin{array}{cc}
\left( A_{F(\lambda)} - \lambda E_H \right)^{-1}, & \lambda\in \Pi_{\lambda_0}\\
\left( A_{F^*(\overline{\lambda})} - \lambda E_H \right)^{-1}, &
\overline{\lambda}\in \Pi_{\lambda_0}
\end{array}
\right.,
\end{equation}
where $F(\lambda)$ is a function from $\mathcal{S}_{a;\lambda_0}(\Pi_{\lambda_0}; \mathcal{N}_{\lambda_0}(A),
\mathcal{N}_{\overline{\lambda_0}}(A))$.
Conversely, an arbitrary function $F(\lambda)\in \mathcal{S}_{a;\lambda_0}(\Pi_{\lambda_0}; \mathcal{N}_{\lambda_0}(A),
\mathcal{N}_{\overline{\lambda_0}}(A))$
defines by relation~(\ref{f11_3_p2_1}) a generalized resolvent
$\mathbf{R}_{s;\lambda}$ of the operator $A$.
Moreover, for different functions from
$\mathcal{S}_{a;\lambda_0}(\Pi_{\lambda_0}; \mathcal{N}_{\lambda_0}(A),
\mathcal{N}_{\overline{\lambda_0}}(A))$
there correspond different generalized resolvents of the operator $A$.
By
$\mathcal{S}_{a; \lambda_0} (\Pi_{\lambda_0}; \mathcal{N}_{\lambda_0}(A),
\mathcal{N}_{ \overline{\lambda_0} }(A)) = \mathcal{S}_{a} (\Pi_{\lambda_0}; \mathcal{N}_{\lambda_0}(A),
\mathcal{N}_{ \overline{\lambda_0} }(A))$
we denote a
set of all operator-valued functions $F(\lambda)\in \mathcal{S}(\Pi_{\lambda_0}; \mathcal{N}_{\lambda_0}(A),
\mathcal{N}_{ \overline{\lambda_0} }(A))$, which are $\lambda_0$-admissible with respect to the operator $A$.
\textit{Throughout this paper we shall use the notations from the survey paper~\cite{cit_1000_Zagorodnyuk}}.

Our main aim here is to characterize those generalized resolvents $\mathbf R_{s;\lambda}$, which are generated by
an (at least one) invertible self-adjoint extension. Such generalized resolvents we shall call
\textbf{generalized I-resolvents}. A generalized I-resolvent will be described by Shtraus's formula~(\ref{f11_3_p2_1}),
where the parameter $F(\lambda)$ satisfies a boundary condition at $0$.
As a by-product, using the generalized Neumann's formulas we describe all invertible extensions of a symmetric operator inside $H$.

{\bf Notations. } As it was already mentioned, in this paper we use notations from~\cite{cit_1000_Zagorodnyuk}.

\section{Invertible extensions.}

Let $A$ be a closed symmetric invertible operator in a Hilbert space $H$.
Let
$z$ from $\mathbb{C}_-$ ($\mathbb{C}_+$) be a fixed number.
Recall that the following formulas (see~\cite[Theorem 3.13]{cit_1000_Zagorodnyuk})
\begin{equation}
\label{f4_1_p2_1}
D(B) = D(A) \dotplus (T-E_H) D(T),
\end{equation}
\begin{equation}
\label{f4_2_p2_1}
B ( f + T\psi - \psi ) = Af + zT\psi - \overline{z}\psi,\qquad f\in D(A),\ \psi\in D(T),
\end{equation}
establish a one-to-one correspondence between all admissible with respect to $A$ isometric
operators $T$, $D(T)\subseteq \mathcal{N}_z(A)$, $R(T)\subseteq \mathcal{N}_{\overline{z}}(A)$,
and all symmetric extensions $B$ of the operator $A$.
We have
\begin{equation}
\label{f4_3_p2_1}
D(T) = \mathcal{N}_z(A)\cap R(B-zE_H),
\end{equation}
\begin{equation}
\label{f4_4_p2_1}
T\subseteq (B-\overline{z}E_H)(B-zE_H)^{-1}.
\end{equation}
Formulas~(\ref{f4_1_p2_1}),(\ref{f4_2_p2_1})
define a one-to-one correspondence between all admissible with respect to $A$ non-expanding
operators $T$, $D(T)\subseteq \mathcal{N}_z(A)$, $R(T)\subseteq \mathcal{N}_{\overline{z}}(A)$,
and all dissipative (respectively accumulative) extensions $B$ of the operator $A$.
Relations~(\ref{f4_3_p2_1}),(\ref{f4_4_p2_1}) hold in this case, as well.

Consider the Cayley transformation of the operator $A$:
\begin{equation}
\label{f3_1_2_p1_1}
U_z = U_z(A) = (A - \overline{z} E_H)(A-zE_H)^{-1} = E_H + (z-\overline{z}) (A-zE_H)^{-1}.
\end{equation}
The operator $B$ may be also determined by the following relations:
\begin{equation}
\label{f2_0_p2_1}
W_z = (B - \overline{z} E_H)(B-zE_H)^{-1} = E_H + (z-\overline{z}) (B-zE_H)^{-1},
\end{equation}
\begin{equation}
\label{f2_0_1_p2_1}
B = (z W_z - \overline{z} E_H) (W_z - E_H)^{-1} = zE_H + (z-\overline{z})(W_z - E_H)^{-1}.
\end{equation}
\begin{equation}
\label{f2_0_2_p2_1}
W_z = U_z(A)\oplus T.
\end{equation}
It is readily checked that
\begin{equation}
\label{k1_1}
\mathcal{M}_\lambda(A) = \mathcal{M}_{ \frac{1}{\lambda} }(A^{-1}),\quad
\mathcal{N}_\lambda(A) = \mathcal{N}_{ \frac{1}{\lambda} }(A^{-1}),\qquad \lambda\in \mathbb{R}_e;
\end{equation}
\begin{equation}
\label{k1_2}
U_\lambda(A) = \frac{ \overline{\lambda} }{ \lambda } U_{ \frac{1}{\lambda} }(A^{-1}),\quad
\lambda\in \mathbb{R}_e.
\end{equation}

\begin{thm}
\label{tk1_1}
Let $A$ be a closed symmetric invertible operator in a Hilbert space $H$, and
$z\in \mathbb{R}_e$ be a fixed point.
Let $T$
be an admissible with respect to $A$ non-expanding
operator with
$D(T)\subseteq \mathcal{N}_z(A)$, $R(T)\subseteq \mathcal{N}_{\overline{z}}(A)$.
The following two conditions are equivalent:

\begin{itemize}
\item[(i)] The operator $B$, defined by~(\ref{f4_1_p2_1}) and (\ref{f4_2_p2_1}), is invertible;

\item[(ii)] The operator $\frac{z}{ \overline{z} } T$ is $\frac{1}{z}$-admissible with respect to $A^{-1}$.

\end{itemize}
\end{thm}
\textbf{Proof.}
From relation~(\ref{f2_0_1_p2_1}) it follows that $B$ is invertible if and only if
the operator $\frac{z}{ \overline{z} } W_z$ has no non-zero fixed elements.
By~(\ref{f2_0_2_p2_1}),(\ref{k1_2}) we see that
$$ \frac{z}{ \overline{z} } W_z = \frac{z}{ \overline{z} } U_z(A)\oplus \frac{z}{ \overline{z} } T =
U_{ \frac{1}{z} }(A^{-1})\oplus \frac{z}{ \overline{z} } T. $$
It remains to apply Theorem~3.6 from~\cite{cit_1000_Zagorodnyuk}.
$\Box$

\begin{rmr}
\label{rk1_1}
By the definition of an admissible operator condition~(ii) of the last theorem is equivalent to the following condition:

\noindent
(iii) the operator $T - \frac{ \overline{z} }{z} X_{ \frac{1}{z} }(A^{-1})$ is invertible,

where $X_{ \lambda }(A^{-1})$ is the forbidden operator for $A^{-1}$.
\end{rmr}

\begin{thm}
\label{tk1_2}
Let $A$ be a closed symmetric invertible operator in a Hilbert space $H$. Suppose that $A$ has finite defect numbers.
Then there exists an invertible self-adjoint operator $\widetilde A\supseteq A$ in a Hilbert space
$\widetilde H\supseteq H$.
\end{thm}
\textbf{Proof.}
Let the operator $A$ from the formulation of the theorem have the deficiency index $(n,m)$, $n,m\in \mathbb{Z}_+$, $n+m\not=0$.
Set $\mathcal{H} = H \oplus H$ and $\mathcal{A} = A \oplus (-A)$. The closed symmetric operator $\mathcal{A}$
has equal defect numbers. Fix an arbitrary number $z\in \mathbb{R}_e$. Let
$\{ f_k \}_{k=1}^{n+m}$ be an orthonormal basis in $\mathcal{N}_z(\mathcal{A})$.
Set
$$ T (\alpha f_1) = \alpha h,\qquad \alpha\in \mathbb{C}, $$
where $h\in \mathcal{N}_{ \overline{z} }(\mathcal{A})$ is an arbitrary element such that $\| h \|_{\mathcal{H}} = 1$, and
$h\not = \frac{ \overline{z} }{z} X_{ \frac{1}{z} }(\mathcal{A}^{-1}) f_1$
if
$f_1\in D(X_{ \frac{1}{z} }(\mathcal{A}^{-1}))$;
$h\not = X_{ z }(\mathcal{A}) f_1$ if
$f_1\in D(X_z(\mathcal{A}))$.
The operator $T$ with the domain $\mathop{\rm span}\nolimits\{ f_1 \}$ is $z$-admissible with respect to $\mathcal{A}$ and
condition~(iii) holds for $\mathcal{A}$.
By Theorem~\ref{tk1_1} the symmetric extension $B$, corresponding to $T$, is invertible. The defect numbers of $B$
are equal to $n+m-1$. If $B$ is not self-adjoint, we can take $B$ instead of $\mathcal{A}$ in the above construction to obtain a closed
symmetric extension with the deficiency index $(n+m-2,n+m-2)$.
Repeating this procedure we shall construct a self-adjoint invertible extension of $\mathcal{A}$.
$\Box$




\section{Generalized resolvents.}

Consider a closed symmetric invertible operator $A$ in a Hilbert space $H$.
Choose and fix an arbitrary point $\lambda_0\in \mathbb{R}_e$.
A function $F(\lambda)\in
\mathcal{S}_{a; \lambda_0} (\Pi_{\lambda_0}; \mathcal{N}_{\lambda_0}(A),
\mathcal{N}_{ \overline{\lambda_0} }(A))$
is said to be  \textbf{$\lambda_0$-I-admissible with respect to the operator $A$},
if the validity of
\begin{equation}
\label{k2_1}
\lim_{\lambda\in\Pi_{\lambda_0}^\varepsilon,\ \lambda\to 0} F(\lambda) \psi = \frac{ \overline{\lambda_0} }{\lambda_0} X_{ \frac{1}{\lambda_0} }(A^{-1})\psi,
\end{equation}
\begin{equation}
\label{k2_2}
\underline{\lim}_{\lambda\in\Pi_{\lambda_0}^\varepsilon,\ \lambda\to 0}
\left[
\frac{1}{ |\lambda| }
(\| \psi \|_H - \| F(\lambda) \psi \|_H)
\right] < +\infty,
\end{equation}
for some $\varepsilon$: $0<\varepsilon <\frac{\pi}{2}$,
implies $\psi = 0$.

A set of all operator-valued functions $F(\lambda)\in
\mathcal{S}_{a; \lambda_0} (\Pi_{\lambda_0}; \mathcal{N}_{\lambda_0}(A),
\mathcal{N}_{ \overline{\lambda_0} }(A))$
which are $\lambda_0$-I-admissible with respect to the operator $A$,
we shall denote by
$$ \mathcal{S}_{I; \lambda_0} (\Pi_{\lambda_0}; \mathcal{N}_{\lambda_0}(A),
\mathcal{N}_{ \overline{\lambda_0} }(A)) = \mathcal{S}_{I} (\Pi_{\lambda_0}; \mathcal{N}_{\lambda_0}(A),
\mathcal{N}_{ \overline{\lambda_0} }(A)). $$

In the case $\overline{R(A)} = H$, we have $D(X_{ \frac{1}{\lambda_0} }(A^{-1})) = \{ 0 \}$ and therefore
$\mathcal{S}_{I;\lambda_0}(\Pi_{\lambda_0}; \mathcal{N}_{\lambda_0}(A),
\mathcal{N}_{ \overline{\lambda_0} }(A)) =
\mathcal{S}_{a;\lambda_0}(\Pi_{\lambda_0}; \mathcal{N}_{\lambda_0}(A),
\mathcal{N}_{ \overline{\lambda_0} }(A))$.

\begin{thm}
\label{kt2_1}
Let $A$ be a closed symmetric invertible operator in a Hilbert space $H$, and $\lambda_0\in \mathbb{R}_e$ be an arbitrary
point.
An arbitrary generalized I-resolvent $\mathbf{R}_{\lambda}$ of the operator $A$ has the following form:
\begin{equation}
\label{k2_3}
\mathbf R_{\lambda} = \left\{ \begin{array}{cc}
\left( A_{F(\lambda)} - \lambda E_H \right)^{-1}, & \lambda\in \Pi_{\lambda_0}\\
\left( A_{F^*(\overline{\lambda})} - \lambda E_H \right)^{-1}, &
\overline{\lambda}\in \Pi_{\lambda_0}
\end{array}
\right.,
\end{equation}
where $F(\lambda)$ is a function from $\mathcal{S}_{I;\lambda_0}(\Pi_{\lambda_0}; \mathcal{N}_{\lambda_0}(A),
\mathcal{N}_{\overline{\lambda_0}}(A))$.
Conversely, an arbitrary function $F(\lambda)\in \mathcal{S}_{I;\lambda_0}(\Pi_{\lambda_0}; \mathcal{N}_{\lambda_0}(A),
\mathcal{N}_{\overline{\lambda_0}}(A))$
defines by relation~(\ref{k2_3}) a generalized I-resolvent
$\mathbf{R}_{\lambda}$ of the operator $A$.
Moreover, for different functions from
$\mathcal{S}_{I;\lambda_0}(\Pi_{\lambda_0}; \mathcal{N}_{\lambda_0}(A),
\mathcal{N}_{\overline{\lambda_0}}(A))$
there correspond different generalized I-resolvents of the operator $A$.
\end{thm}
\textbf{Proof.}
Let $A$ be a closed symmetric invertible operator in a Hilbert space $H$, $\lambda_0\in \mathbb{R}_e$.

Let us prove the first statement of the theorem.
Let $\mathbf{R}_{\lambda}$ be an arbitrary generalized I-resolvent of the operator $A$.
The generalized I-resolvent $\mathbf{R}_{\lambda}$ is generated by a self-adjoint invertible operator
$\widetilde A\supseteq A$ in a Hilbert space $\widetilde H\supseteq H$.
Repeating the arguments from the proof of Theorem~3.34 in~\cite{cit_1000_Zagorodnyuk}
we obtain a representation~(\ref{k2_3}), where $F(\lambda)$
is a function from $\mathcal{S}_{a;\lambda_0}(\Pi_{\lambda_0}; \mathcal{N}_{\lambda_0}(A),
\mathcal{N}_{\overline{\lambda_0}}(A))$.
It remains to check that $F(\lambda)\in \mathcal{S}_{I;\lambda_0}(\Pi_{\lambda_0}; \mathcal{N}_{\lambda_0}(A),
\mathcal{N}_{\overline{\lambda_0}}(A))$.

The self-adjoint extension $\widetilde A^{-1}$ of the operator $A^{-1}$ in a Hilbert space $\widetilde H\supseteq H$
generates the following objects, see Subsection~3.6 in~\cite{cit_1000_Zagorodnyuk}:
the operator-valued function $\mathfrak{B}_\lambda (A^{-1},\widetilde A^{-1})$, $\lambda\in \mathbb{R}_e$;
the operator $\mathfrak{B}_\infty(A^{-1},\widetilde A^{-1})$;
the operator-valued function $\mathfrak{F}(\lambda; \frac{1}{\lambda_0}, A^{-1}, \widetilde A^{-1})$, $\lambda\in \Pi_{ \frac{1}{\lambda_0} }$;
the operator $\Phi_\infty(\frac{1}{\lambda_0}; A^{-1},\widetilde A^{-1})$.

On the other hand, the operator $A$ can be identified with the operator
\begin{equation}
\label{k2_4}
\mathcal{A} = A\oplus A_e,
\end{equation}
with $A_e = o_H$, in the Hilbert space
\begin{equation}
\label{k2_5}
\widetilde{H} = H\oplus H_e,
\end{equation}
where $H_e := \widetilde{H} \ominus H$, see Subsection~3.5 in~\cite{cit_1000_Zagorodnyuk}.
Then
\begin{equation}
\label{k2_6}
\mathcal{A}^{-1} = A^{-1}\oplus A_e^{-1} = A^{-1}\oplus A_e.
\end{equation}
By the generalized Neumann's formulas, for $\widetilde A\supseteq \mathcal{A}$ there corresponds an isometric operator
$T$, $D(T) = \mathcal{N}_{\lambda_0}(\mathcal{A})$, $R(T) = \mathcal{N}_{ \overline{\lambda_0} }(\mathcal{A})$,
which is $\lambda_0$-admissible with respect to $\mathcal{A}$. Moreover, since the operator $\widetilde A$ is invertible,
then by Theorem~\ref{tk1_1} we conclude that
the operator $\frac{\lambda_0}{ \overline{\lambda_0} } T$ is
$\frac{1}{\lambda_0}$-admissible with respect to $\mathcal{A}^{-1}$.

Applying Theorem~3.16 in~\cite{cit_1000_Zagorodnyuk} for the operator $\mathcal{A}^{-1}$, with $z=\frac{1}{\lambda_0}$, we
obtain that the operator $\Phi(\frac{1}{\lambda_0};\mathcal{A}^{-1},\frac{\lambda_0}{ \overline{\lambda_0} } T)$
is $\frac{1}{\lambda_0}$-admissible with respect to $A^{-1}$.

By the generalized Neumann's formulas, for $\widetilde A^{-1}\supseteq \mathcal{A}^{-1}$ there corresponds an isometric operator
$V$, $D(V) = \mathcal{N}_{\frac{1}{\lambda_0}}(\mathcal{A}^{-1}) =
\mathcal{N}_{\lambda_0}(\mathcal{A})$,
$R(V) = \mathcal{N}_{ \frac{1}{\overline{\lambda_0}} }(\mathcal{A}^{-1}) = \mathcal{N}_{ \overline{\lambda_0} }(\mathcal{A})$,
which is $\frac{1}{\lambda_0}$-admissible with respect to $\mathcal{A}^{-1}$.
Consider the following operator $\mathfrak{B}(\frac{1}{\lambda_0};\mathcal{A}^{-1},V)$, see~(3.38) in~\cite{cit_1000_Zagorodnyuk}:
$$\mathfrak{B}(\frac{1}{\lambda_0};\mathcal{A}^{-1},V) h
= P^{\widetilde H}_H (\mathcal{A}^{-1})_V h = P^{\widetilde H}_H \widetilde A^{-1} h,\qquad h\in D((\mathcal{A}^{-1})_V) \cap H = D(\widetilde A^{-1})\cap H. $$
Comparing this definition with the definition of the operator $\mathfrak{B}_\infty(A^{-1},\widetilde A^{-1})$
we conclude that
\begin{equation}
\label{k2_7}
\mathfrak{B}(\frac{1}{\lambda_0};\mathcal{A}^{-1},V) = \mathfrak{B}_\infty(A^{-1},\widetilde A^{-1}).
\end{equation}
Applying Theorem~3.20 in~\cite{cit_1000_Zagorodnyuk} for the operator $\widetilde{A}^{-1}$, with $z=\frac{1}{\lambda_0}$, we get:
$$  \mathfrak{B}(\frac{1}{\lambda_0};\mathcal{A}^{-1},V) = (A^{-1})_{\Phi(\frac{1}{\lambda_0};\mathcal{A}^{-1}, V), \frac{1}{\lambda_0}}. $$
On the other hand, by the definition of $\Phi_\infty(\frac{1}{\lambda_0}; A^{-1},\widetilde A^{-1})$ we have:
$$ \mathfrak{B}_\infty(A^{-1},\widetilde A^{-1}) = (A^{-1})_{\Phi_\infty(\frac{1}{\lambda_0}; A^{-1},\widetilde A^{-1}), \frac{1}{\lambda_0}}. $$
Then
\begin{equation}
\label{k2_8}
\Phi_\infty(\frac{1}{\lambda_0}; A^{-1},\widetilde A^{-1}) = \Phi(\frac{1}{\lambda_0};\mathcal{A}^{-1}, V).
\end{equation}

By the generalized Neumann's formulas~(\ref{f4_3_p2_1}),(\ref{f4_4_p2_1}) we get:
$$ D(T) = \mathcal{N}_{\lambda_0}(\mathcal{A})\cap R(\widetilde A - \lambda_0 E_{\widetilde H}), $$
$$ T\subseteq (\widetilde A - \overline{\lambda_0} E_{\widetilde H}) (\widetilde A - \lambda_0 E_{\widetilde H})^{-1}; $$
$$ D(V) = \mathcal{N}_{\frac{1}{\lambda_0}}(\mathcal{A}^{-1})\cap R(\widetilde A^{-1} - \frac{1}{\lambda_0} E_{\widetilde H}) =
\mathcal{N}_{\lambda_0}(\mathcal{A})\cap R(\widetilde A - \lambda_0 E_{\widetilde H}), $$
$$ V\subseteq (\widetilde A^{-1} - \frac{1}{\overline{\lambda_0}} E_{\widetilde H}) (\widetilde A^{-1} - \frac{1}{\lambda_0} E_{\widetilde H})^{-1} $$
$$ = \frac{\lambda_0}{\overline{\lambda_0}}
(\widetilde A - \overline{\lambda_0} E_{\widetilde H}) (\widetilde A - \lambda_0 E_{\widetilde H})^{-1}. $$
Therefore
$$ V = \frac{\lambda_0}{\overline{\lambda_0}} T. $$
By~(\ref{k2_8}) we get:
\begin{equation}
\label{k2_9}
\Phi_\infty(\frac{1}{\lambda_0}; A^{-1},\widetilde A^{-1}) = \Phi(\frac{1}{\lambda_0};\mathcal{A}^{-1}, \frac{\lambda_0}{\overline{\lambda_0}} T).
\end{equation}
Therefore the operator $\Phi_\infty(\frac{1}{\lambda_0}; A^{-1},\widetilde A^{-1})$
is $\frac{1}{\lambda_0}$-admissible with respect to $A^{-1}$.

The self-adjoint operator $\widetilde A\supseteq A$ in a Hilbert space $\widetilde H\supseteq H$
generates the operator-valued function $\mathfrak{B}_\lambda (A,\widetilde A)$, $\lambda\in \mathbb{R}_e$,
and the operator-valued function $\mathfrak{F}(\lambda; \lambda_0, A, \widetilde A)$, $\lambda\in \Pi_{\lambda_0}$
(see~\cite[Subsection 3.6]{cit_1000_Zagorodnyuk}).
Denote (see~\cite[p. 222]{cit_1000_Zagorodnyuk})
\begin{equation}
\label{f6_3_p2_1}
\widetilde{\mathfrak{L}}_\lambda (A,\widetilde A) =
\left\{
h\in D(\widetilde A):\ (\widetilde A - \lambda E_{\widetilde H})h \in H
\right\},\qquad \lambda\in \mathbb{C};
\end{equation}
$$
\mathfrak{L}_\lambda (A,\widetilde A) = P^{\widetilde H}_H   \widetilde{\mathfrak{L}}_\lambda,\qquad \lambda\in \mathbb{C}.
$$
We shall also need sets $\widetilde{\mathfrak{L}}_\lambda (A^{-1},\widetilde A^{-1})$, $\mathfrak{L}_\lambda (A^{-1},\widetilde A^{-1})$, $\lambda\in \mathbb{C}$,
constructed by~(\ref{f6_3_p2_1}) with the operators $A^{-1}$,$\widetilde A^{-1}$ instead of $A$,$\widetilde A$.

Choose an arbitrary element $h\in \widetilde{\mathfrak{L}}_\lambda (A,\widetilde A)$, $\lambda\in \mathbb{R}_e$. Then $\widetilde A h\in R(\widetilde A)$, and
$$ (\widetilde A^{-1} - \frac{1}{\lambda} E_{\widetilde H}) \widetilde A h = -\frac{1}{\lambda}
(\widetilde A - \lambda E_{\widetilde H}) h\in H. $$
Therefore
$\widetilde A h\in \widetilde{\mathfrak{L}}_{\frac{1}{\lambda}} (A^{-1},\widetilde A^{-1})$, and
$$ \widetilde A \widetilde{\mathfrak{L}}_\lambda (A,\widetilde A) \subseteq \widetilde{\mathfrak{L}}_{\frac{1}{\lambda}} (A^{-1},\widetilde A^{-1}),\qquad
\lambda\in \mathbb{R}_e. $$
In order to obtain the equality:
\begin{equation}
\label{k2_10}
\widetilde A \widetilde{\mathfrak{L}}_\lambda (A,\widetilde A) = \widetilde{\mathfrak{L}}_{\frac{1}{\lambda}} (A^{-1},\widetilde A^{-1}),\qquad
\lambda\in \mathbb{R}_e,
\end{equation}
it remains to apply the proved inclusion for the operators $A^{-1}$, $\widetilde A^{-1}$ instead of $A$, $\widetilde A$,
and with $\frac{1}{\lambda}$ instead of $\lambda$.

The operator $P^{\widetilde H}_H|_{\widetilde{\mathfrak{L}}_\lambda (A,\widetilde A)}$, $\lambda\in \mathbb{R}_e$,
is invertible, see considerations below~(3.48)
in~\cite{cit_1000_Zagorodnyuk}.
By the definition of
$\mathfrak{B}_\lambda (A,\widetilde A)$ we may write:
\begin{equation}
\label{k2_11}
\mathfrak{B}_\lambda (A,\widetilde A) = P^{\widetilde H}_H \widetilde A
\left(
P^{\widetilde H}_H|_{\widetilde{\mathfrak{L}}_\lambda (A,\widetilde A)}
\right)^{-1},\qquad \lambda\in \mathbb{R}_e.
\end{equation}
Applying this representation for the operators $A^{-1}$, $\widetilde A^{-1}$ instead of $A$, $\widetilde A$,
and with $\frac{1}{\lambda}$ instead of $\lambda$ we get:
$$
\mathfrak{B}_{ \frac{1}{\lambda} } (A^{-1},\widetilde A^{-1}) = P^{\widetilde H}_H \widetilde A^{-1}
\left(
P^{\widetilde H}_H|_{\widetilde{\mathfrak{L}}_{ \frac{1}{\lambda} } (A^{-1},\widetilde A^{-1})}
\right)^{-1},\qquad \lambda\in \mathbb{R}_e.
$$
Observe that
$$ R(\mathfrak{B}_\lambda(A,\widetilde A)) = P^{\widetilde H}_H \widetilde A \widetilde{\mathfrak{L}}_\lambda (A,\widetilde A) =
\mathfrak{L}_{\frac{1}{\lambda}} (A^{-1},\widetilde A^{-1}),\qquad \lambda\in \mathbb{R}_e, $$
and
$$ \mathfrak{B}_{ \frac{1}{\lambda} } (A^{-1},\widetilde A^{-1}) \mathfrak{B}_\lambda(A,\widetilde A) g = g,\qquad
g\in \mathfrak{L}_\lambda(A,\widetilde A),\ \lambda\in \mathbb{R}_e. $$
Therefore
\begin{equation}
\label{k2_12}
\mathfrak{B}_\lambda (A,\widetilde A)^{-1} =
\mathfrak{B}_{\frac{1}{\lambda}} (A^{-1},\widetilde A^{-1}),\qquad \lambda\in \mathbb{R}_e.
\end{equation}
By the definition of the operator-valued function $\mathfrak{F}(\lambda)$ and by~(\ref{k2_12}) we may write:
$$
\mathfrak{F}(\lambda; \lambda_0, A, \widetilde A) =
\left. (\mathfrak{B}_\lambda(A, \widetilde A) - \overline{\lambda_0}E_H)
(\mathfrak{B}_\lambda(A, \widetilde A) - \lambda_0 E_H)^{-1} \right|_{\mathcal{N}_{\lambda_0}(A)},\qquad \lambda\in \Pi_{\lambda_0};
$$
and, also for $\lambda\in \Pi_{ \lambda_0 }$, we have:
$$
\mathfrak{F}(\frac{1}{\lambda}; \frac{1}{\lambda_0}, A^{-1}, \widetilde A^{-1}) $$
$$ =
\left. (\mathfrak{B}_{ \frac{1}{\lambda} } (A^{-1},\widetilde A^{-1}) - \frac{1}{\overline{\lambda_0}} E_H)
(\mathfrak{B}_{ \frac{1}{\lambda} } (A^{-1},\widetilde A^{-1}) - \frac{1}{\lambda_0} E_H)^{-1} \right|_{\mathcal{N}_{\lambda_0}(A)}
$$
$$ =
\left.
\frac{\lambda_0}{ \overline{\lambda_0} }
(\mathfrak{B}_{ \lambda } (A,\widetilde A) - \overline{\lambda_0} E_H)
(\mathfrak{B}_{ \lambda } (A,\widetilde A) - \lambda_0 E_H)^{-1} \right|_{\mathcal{N}_{\lambda_0}(A)}.
$$

Therefore
\begin{equation}
\label{k2_13}
\mathfrak{F}(\frac{1}{\lambda}; \frac{1}{\lambda_0}, A^{-1}, \widetilde A^{-1}) =
\frac{\lambda_0}{ \overline{\lambda_0} }
\mathfrak{F}(\lambda; \lambda_0, A, \widetilde A),\qquad \lambda\in \Pi_{\lambda_0}.
\end{equation}

Suppose that relations~(\ref{k2_1}),(\ref{k2_2}) with $F(\lambda) = \mathfrak{F}(\lambda; \lambda_0, A, \widetilde A)$ hold
for some $\varepsilon$: $0<\varepsilon <\frac{\pi}{2}$.
Then using the change of a variable $y=\frac{1}{\lambda}$ we get:
\begin{equation}
\label{k2_14}
\lim_{y\in\Pi_{ \frac{1}{\lambda_0} }^\varepsilon,\ y\to \infty} \mathfrak{F}(\frac{1}{y}; \lambda_0, A, \widetilde A) \psi =
\frac{ \overline{\lambda_0} }{\lambda_0} X_{ \frac{1}{\lambda_0} }(A^{-1})\psi,
\end{equation}
\begin{equation}
\label{k2_15}
\underline{\lim}_{y\in\Pi_{  \frac{1}{\lambda_0}  }^\varepsilon,\ y\to \infty}
\left[
|y|
(\| \psi \|_H - \| \mathfrak{F}(\frac{1}{y}; \lambda_0, A, \widetilde A) \psi \|_H)
\right] < +\infty.
\end{equation}
By~(\ref{k2_13}) we obtain that
$$
\mathfrak{F}(\frac{1}{y}; \lambda_0, A, \widetilde A)
=
\frac{ \overline{\lambda_0} }{\lambda_0}
\mathfrak{F}(y; \frac{1}{\lambda_0}, A^{-1}, \widetilde A^{-1}),\qquad y\in \Pi_{\frac{1}{\lambda_0}}.
$$
Substituting this expression in relations~(\ref{k2_14}),(\ref{k2_15}) we get:
\begin{equation}
\label{k2_16}
\lim_{y\in\Pi_{ \frac{1}{\lambda_0} }^\varepsilon,\ y\to \infty}
\mathfrak{F}(y; \frac{1}{\lambda_0}, A^{-1}, \widetilde A^{-1}) \psi =
 X_{ \frac{1}{\lambda_0} }(A^{-1})\psi,
\end{equation}
\begin{equation}
\label{k2_17}
\underline{\lim}_{y\in\Pi_{  \frac{1}{\lambda_0}  }^\varepsilon,\ y\to \infty}
\left[
|y|
(\| \psi \|_H - \|
\mathfrak{F}(y; \frac{1}{\lambda_0}, A^{-1}, \widetilde A^{-1}) \psi \|_H)
\right] < +\infty.
\end{equation}

By Theorem~3.32 in~\cite{cit_1000_Zagorodnyuk} we conclude that
$\psi\in D(\Phi_\infty(\frac{1}{\lambda_0};A^{-1},\widetilde A^{-1}))$, and
$$ \Phi_\infty(\frac{1}{\lambda_0};A^{-1},\widetilde A^{-1}) \psi = X_{ \frac{1}{\lambda_0} }(A^{-1})\psi. $$
Since $\Phi_\infty(\frac{1}{\lambda_0};A^{-1},\widetilde A^{-1})$ is
$\frac{1}{\lambda_0}$-admissible with respect to $A^{-1}$, we obtain that $\psi = 0$.
Consequently, $F(\lambda)\in \mathcal{S}_{I;\lambda_0}(\Pi_{\lambda_0}; \mathcal{N}_{\lambda_0}(A),
\mathcal{N}_{\overline{\lambda_0}}(A))$.

Let us check the second statement of the theorem.
Let $F(\lambda)$ be an arbitrary function from $\mathcal{S}_{I;\lambda_0}(\Pi_{\lambda_0}; \mathcal{N}_{\lambda_0}(A),
\mathcal{N}_{\overline{\lambda_0}}(A))$.
We shall use the following lemma.

\begin{lem}
\label{lk2_1}
Let $n,m\in \mathbb{Z}_+\cup\{ \infty \}$: $n+m\not= 0$.
There exists a closed symmetric invertible operator $A$ in a Hilbert space $H$, $\overline{D(A)} = H$,
$\overline{R(A)} = H$, which has the deficiency index $(n,m)$.
\end{lem}
\textbf{Proof.}
Let $H_0$ be an arbitrary Hilbert space, and $\mathfrak{A} = \{ f_k \}_{k=0}^\infty$ be an orthonormal basis in $H_0$.
Consider the following operator (unilateral shift):
$$ V_0 h = \sum_{k=0}^\infty \alpha_k f_{k+1},\qquad h=\sum_{k=0}^\infty \alpha_k f_k\in H_0,\ \alpha_k\in \mathbb{C}, $$
with $D(V_0) = H_0$. The operator $V_0$ is closed, isometric, and its deficiency index is $(0,1)$.
The condition $V_0 g = \pm g$, for an element $g\in H_0$, implies $g=0$.
Consequently, the inverse Cayley transformation:
$$ A_0 = i ( V_0 + E_{H_0} ) ( V_0 - E_{H_0} )^{-1}, $$
is a closed symmetric invertible operator in $H_0$, with the deficiency index $(0,1)$.
If $h\in H_0$ and $h\perp ( V_0 \pm E_{H_0} ) H_0$, then $V_0^* h = \mp h$.
The condition $V_0^* h = V_0^{-1} P^{H_0}_{R(V_0)} h = \mp h$, for an element $h\in H_0$, implies $h=0$.
Therefore, $\overline{D(A_0)} = H_0$, $\overline{R(A_0)} = H_0$.
Set
$$ H^l = \bigoplus\limits_{j=0}^l H_0,\quad W_l = \bigoplus\limits_{j=0}^l V_0,\quad l\in \mathbb{Z}_+\cup\{ \infty \}. $$
$W_l$ is a closed isometric operator in $H^l$.
The deficiency index of $W_l$ is equal to $(0,l+1)$. If $W_l h = \pm h$, or $W_l^* h = \pm h$, then $h=0$.
Then
$$ A_l = i ( W_l + E_{H^l} ) ( W_l - E_{H^l} )^{-1}, $$
is a closed symmetric invertible operator in $H^l$, with the deficiency index $(0,l+1)$, $l\in \mathbb{Z}_+$.
Moreover, we have $\overline{D(A_l)} = H^l$, $\overline{R(A_l)} = H^l$.
Observe that the operator $-A_l$ has the deficiency index $(l+1,0)$, $l\in \mathbb{Z}_+$.
If $m>0,n>0$, we set
$$ H = H^{m-1} \oplus H^{n-1},\quad A = (-A_{m-1}) \oplus A_{n-1}. $$
$\Box$

Let us return to the proof of the theorem.
If the operator $A$ is self-adjoint, then the set $\mathcal{S}_{I;\lambda_0}(\Pi_{\lambda_0}; \mathcal{N}_{\lambda_0}(A),
\mathcal{N}_{\overline{\lambda_0}}(A))$ consists of a unique function $F(\lambda)=o_H$. Of course, this function generates
the resolvent of $A$ by~(\ref{k2_3}).
Thus, we can assume that $A$ is not self-adjoint.

We shall use the scheme of the proof of the corresponding statement in Theorem~3.34 in~\cite{cit_1000_Zagorodnyuk}.
By~Lemma~\ref{lk2_1} there exists a closed symmetric invertible operator $A_1$ in a Hilbert space $H_1$,
$\overline{D(A_1)} = H_1$, $\overline{R(A_1)} = H_1$, which has the same defect numbers as $A$.
Let $U$ and $W$ be arbitrary isometric operators, which map respectively $\mathcal{N}_{\lambda_0}(A_1)$ on $\mathcal{N}_{\lambda_0}(A)$,
and $\mathcal{N}_{\overline{\lambda_0}}(A_1)$ on $\mathcal{N}_{\overline{\lambda_0}}(A)$.
Set $F_1(\lambda) = W^{-1} F(\lambda) U$, $\lambda\in\Pi_{\lambda_0}$.
Since $F_1(\lambda)$ belongs to
$\mathcal{S}(\Pi_{\lambda_0}; \mathcal{N}_{\lambda_0}(A_1),
\mathcal{N}_{\overline{\lambda_0}}(A_1))$, by Shtraus's formula it generates a generalized resolvent $\mathbf{R}_\lambda(A_1)$ of $A_1$.

Let us check that $\mathbf{R}_\lambda(A_1)$ is generated by a self-adjoint
\textit{invertible}
operator $\widetilde A_1\supseteq A_1$ in a Hilbert space $\widetilde H_1\supseteq H_1$.
Suppose that $\mathbf{R}_\lambda(A_1)$ is generated by a self-adjoint
operator $\widehat A_1\supseteq A_1$ in a Hilbert space $\widehat H_1\supseteq H_1$.
Suppose that $U_i(\widehat A_1) h = -h$, for an element $h\in \widehat H_1$. Then
$$ (h, U_i(\widehat A_1) g)_{\widehat H_1} = - (U_i(\widehat A_1) h, U_i(\widehat A_1) g)_{\widehat H_1} = - ( h, g)_{\widehat H_1},\qquad
g\in\widehat H_1, $$
$$ 0 = (h, (U_i(\widehat A_1) + E_{\widehat H_1}) g)_{\widehat H_1},\qquad g\in \widehat H_1. $$
In particular, $h$ is orthogonal to $(U_i(A_1) + E_{H_1}) D(U_i(A_1)) = R(A_1)$. Then $h\in \widehat H_1\ominus H_1$.
Set
$$ \widehat H_0 = \{ h\in \widehat H_1:\ U_i(\widehat A_1) h = - h \}. $$
Observe that $\widehat H_0$ is a subspace of $\widehat H_1\ominus H_1$.
Then
$$ \widehat H_1 = H_1\oplus (\widehat H_1\ominus H_1)
=
H_1\oplus \left(
(\widehat H_1\ominus H_1) \ominus \widehat H_0
\right)
\oplus \widehat H_0 = \widetilde H_1 \oplus \widehat H_0, $$
where $\widetilde H_1 := H_1\oplus \left(
(\widehat H_1\ominus H_1) \ominus \widehat H_0
\right)$.

Notice that $U_i(\widehat A_1) \widehat H_0 = \widehat H_0$ and $U_i(\widehat A_1) \widetilde H_1 = \widetilde H_1$.
Set $W_1 = U_i(\widehat A_1)|_{\widetilde H_1}$. There are no non-zero elements $g\in\widetilde H_1$ such that $W_1 g = - g$.
Then the inverse Cayley transformation
$$ \widetilde A_1 := i(W_1 + E_{\widetilde H_1}) (W_1 - E_{\widetilde H_1})^{-1}, $$
is invertible.
Since $\widetilde A_1\subseteq \widehat A_1$, then
$$ (\widetilde A_1 - \lambda E_{\widetilde H_1})^{-1} \subseteq (\widehat A_1 - \lambda E_{\widehat H_1})^{-1}. $$
Therefore $\widetilde A_1$ generates $\mathbf{R}_\lambda(A_1)$.

Set
$$ H_e := \widetilde H_1\ominus H_1. $$

Starting from the same formula, we repeat the rest of the arguments in the proof of Theorem~3.34 in~\cite{cit_1000_Zagorodnyuk}.
In what follows, \textit{ we shall use notations and constructions from this proof without additional references}.
We shall obtain
a self-adjoint operator $\widetilde A\supseteq A$ in a Hilbert space $\widetilde H\supseteq H$, which
generates a generalized resolvent $\mathbf{R}_\lambda$ of $A$. This generalized resolvent is related to $F(\lambda)$
by~(\ref{k2_3}). It remains to check that the operator $\widetilde A$ is invertible.

Since the operator $\widetilde A_1$ is invertible, then by Theorem~\ref{tk1_1}
we obtain that
the operator $\frac{\lambda_0}{ \overline{\lambda_0} } T$ is $\frac{1}{\lambda_0}$-admissible with respect to $\mathcal{A}_1^{-1}$.
By Theorem~3.16 in~\cite{cit_1000_Zagorodnyuk} we conclude that
the operators $\Phi(\frac{1}{\lambda_0}; \mathcal{A}_1^{-1}, \frac{\lambda_0}{ \overline{\lambda_0} } T)$
and $\frac{\lambda_0}{ \overline{\lambda_0} } T_{22}$ are
$\frac{1}{\lambda_0}$-admissible with respect to
$A_1^{-1}$ and $A_e^{-1}=A_e$, respectively.

Comparing the domains of $\Phi(\frac{1}{\lambda_0}; \mathcal{A}_1^{-1}, \frac{\lambda_0}{ \overline{\lambda_0} } T)$ and
$\Phi(\frac{1}{\lambda_0}; \mathcal{A}^{-1}, \frac{\lambda_0}{ \overline{\lambda_0} } V)$ we conclude that
\begin{equation}
\label{k2_17_1}
D(\Phi(\frac{1}{\lambda_0}; \mathcal{A}^{-1}, \frac{\lambda_0}{ \overline{\lambda_0} } V)) = U D(\Phi(\frac{1}{\lambda_0}; \mathcal{A}_1^{-1}, \frac{\lambda_0}{ \overline{\lambda_0} } T)).
\end{equation}
Using Remark~3.15 and formula~(3.28) in~\cite[p. 218]{cit_1000_Zagorodnyuk} for
$\Phi(\frac{1}{\lambda_0}; \mathcal{A}^{-1}, \frac{\lambda_0}{ \overline{\lambda_0} } V)$ and
$\Phi(\frac{1}{\lambda_0}; \mathcal{A}_1^{-1}, \frac{\lambda_0}{ \overline{\lambda_0} } T)$, we get:
\begin{equation}
\label{k2_18}
\Phi(\frac{1}{\lambda_0}; \mathcal{A}^{-1}, \frac{\lambda_0}{ \overline{\lambda_0} } V)
=
W
\Phi(\frac{1}{\lambda_0}; \mathcal{A}_1^{-1}, \frac{\lambda_0}{ \overline{\lambda_0} } T)
U^{-1}.
\end{equation}

We can apply the arguments in the proof of the already proved first statement of the theorem for
the operator $A:=A_1$; the point $\lambda_0$; the generalized $I$-resolvent $\mathbf{R}_\lambda := \mathbf{R}_\lambda(A_1)$ of $A_1$,
which is generated by the self-adjoint invertible operator $\widetilde A := \widetilde A_1$ in $\widetilde H_1\supseteq H_1$.

\begin{equation}
\label{k2_20}
\Phi_\infty(\frac{1}{\lambda_0}; A_1^{-1},\widetilde A_1^{-1}) = \Phi(\frac{1}{\lambda_0};\mathcal{A}_1^{-1}, \frac{\lambda_0}{\overline{\lambda_0}} T).
\end{equation}

\begin{equation}
\label{k2_13_new}
\mathfrak{F}(\frac{1}{\lambda}; \frac{1}{\lambda_0}, A_1^{-1}, \widetilde A_1^{-1}) =
\frac{\lambda_0}{ \overline{\lambda_0} }
\mathfrak{F}(\lambda; \lambda_0, A_1, \widetilde A_1),\qquad \lambda\in \Pi_{\lambda_0}.
\end{equation}

Then
\begin{equation}
\label{k2_22_1}
\mathfrak{F}(y; \frac{1}{\lambda_0}, A_1^{-1}, \widetilde A_1^{-1}) =
\frac{\lambda_0}{ \overline{\lambda_0} }
\mathfrak{F}(\frac{1}{y}; \lambda_0, A_1, \widetilde A_1),\qquad y\in \Pi_{ \frac{1}{\lambda_0} }.
\end{equation}

By Theorem~3.32 in~\cite{cit_1000_Zagorodnyuk} we have:
$$ D(\Phi_\infty(\frac{1}{\lambda_0}; A_1^{-1}, \widetilde A_1^{-1})) =
\{ \psi\in \mathcal{N}_{\frac{1}{\lambda_0}}(A_1^{-1}): $$
\begin{equation}
\label{k2_21}
\left.
\underline{\lim}_{ \lambda\in\Pi_{ \frac{1}{\lambda_0} }^\varepsilon, \lambda\to\infty }
\left[
|\lambda| ( \| \psi \|_{H_1} - \| \mathfrak{F}(\lambda; \frac{1}{\lambda_0}, A_1^{-1}, \widetilde A_1^{-1}) \psi \|_{H_1} )
\right] < +\infty
\right\},
\end{equation}
\begin{equation}
\label{k2_22}
\Phi_\infty(\frac{1}{\lambda_0}; A_1^{-1},\widetilde A_1^{-1}) \psi =
\lim_{ \lambda\in\Pi_{ \frac{1}{\lambda_0} }^\varepsilon, \lambda\to\infty }
\mathfrak{F}(\lambda; \frac{1}{\lambda_0}, A_1^{-1}, \widetilde A_1^{-1}) \psi,\ \psi\in D(\Phi_\infty(\frac{1}{\lambda_0}; A_1^{-1},\widetilde A_1^{-1})),
\end{equation}
where $0 < \varepsilon < \frac{\pi}{2}$.

Using~(\ref{k2_17_1}),(\ref{k2_20}),(\ref{k2_21}),(\ref{k2_22}),(\ref{k2_22_1}) and the change of a variable: $y=\frac{1}{\lambda}$, we obtain that
$$ D(\Phi(\frac{1}{\lambda_0}; \mathcal{A}^{-1}, \frac{\lambda_0}{ \overline{\lambda_0} } V)) =
\{ \psi\in \mathcal{N}_{\lambda_0}(A): $$
\begin{equation}
\label{k2_23}
\left.
\underline{\lim}_{ y\in\Pi_{\lambda_0}^\varepsilon, y\to 0 }
\left[
\frac{1}{|y|} ( \| \psi \|_{H} - \| F(y) \psi \|_{H} )
\right] < +\infty
\right\};
\end{equation}
\begin{equation}
\label{k2_24}
\Phi(\frac{1}{\lambda_0}; \mathcal{A}^{-1}, \frac{\lambda_0}{ \overline{\lambda_0} } V) \psi = \frac{\lambda_0}{ \overline{\lambda_0} }
\lim_{ y\in\Pi_{\lambda_0}^\varepsilon, y\to 0 } F(y) \psi,\qquad \psi\in D(\Phi(\frac{1}{\lambda_0}; \mathcal{A}^{-1}, \frac{\lambda_0}{ \overline{\lambda_0} } V)).
\end{equation}

Suppose that there exists an element $\psi\in
D(\Phi(\frac{1}{\lambda_0}; \mathcal{A}^{-1}, \frac{\lambda_0}{ \overline{\lambda_0} } V))
\cap
X_{ \frac{1}{\lambda_0} }(A^{-1})$ such that the following equality holds:
$$ \Phi(\frac{1}{\lambda_0}; \mathcal{A}^{-1}, \frac{\lambda_0}{ \overline{\lambda_0} } V) \psi = X_{ \frac{1}{\lambda_0} }(A^{-1}) \psi. $$
By~(\ref{k2_23}),(\ref{k2_24}) this means that $\psi\in \mathcal{N}_{\lambda_0}(A)$ and
$$ \underline{\lim}_{ y\in\Pi_{\lambda_0}^\varepsilon, y\to 0 }
\left[
\frac{1}{|y|} ( \| \psi \|_{H} - \| F(y) \psi \|_{H} )
\right] < +\infty, $$
$$ \frac{\lambda_0}{ \overline{\lambda_0} }
\lim_{ y\in\Pi_{\lambda_0}^\varepsilon, y\to 0 } F(y) \psi = X_{ \frac{1}{\lambda_0} }(A^{-1}) \psi. $$
Since $F(\lambda)$ is $\lambda_0$-I-admissible with respect to the operator $A$, we get $\psi = 0$.

\noindent
This means that $\Phi(\frac{1}{\lambda_0}; \mathcal{A}^{-1}, \frac{\lambda_0}{ \overline{\lambda_0} } V)$ is $\frac{1}{\lambda_0}$-admissible with respect to $A^{-1}$.
Since $\frac{\lambda_0}{ \overline{\lambda_0} } T_{22}$ is
$\frac{1}{\lambda_0}$-admissible with respect to
$A_e$, then by Theorem~3.16 in~\cite{cit_1000_Zagorodnyuk} we obtain that
the operator $\frac{\lambda_0}{ \overline{\lambda_0} } V$ is $\frac{1}{\lambda_0}$-admissible with respect to
$\mathcal{A}^{-1}$.
By Theorem~\ref{tk1_1} we conclude that the operator $\mathcal{A}_V = \widetilde A$ is invertible.

The last statement of the theorem follows directly from Shtraus's formula.
$\Box$

\begin{center}
{\large\bf Invertible extensions of symmetric operators and the corresponding generalized resolvents.}
\end{center}
\begin{center}
{\bf S.M. Zagorodnyuk}
\end{center}

In this paper we study invertible extensions of a symmetric operator in a Hilbert space $H$.
All such extensions are characterized by a parameter in the generalized Neumann's formulas.
Generalized resolvents, which are generated by the invertible extensions, are extracted by
a boundary condition among all generalized resolvents in the Shtraus formula.

}
\end{document}